%
%
%
%
%
%
%
%
%
\magnification=\magstephalf      
%
%
\vsize=7.5truein                 
\hsize=5.2truein                 
\newskip\stdskip                 
\stdskip=6pt plus3pt minus3pt    
\medskipamount=\stdskip          
\parindent=0pt                   
\parskip=\stdskip                
\abovedisplayskip=\stdskip       
\belowdisplayskip=\stdskip       
\mathsurround=0.75pt             
\overfullrule=0pt                
%
%
\def\ppar{\par\goodbreak\vskip 8pt plus 4pt minus 4pt}     
%
%
\def\stdspace{\hskip 0.75em plus 0.15em\ignorespaces}
\let\qua\stdspace 
%
%
%
%
%
%
%
\def\hexnumber#1{\ifcase#1 0\or 1\or 2\or 3\or 4\or 5\or 6\or 7\or 8\or
 9\or A\or B\or C\or D\or E\or F\fi}
%
%
\font\thirtnmsa=msam10 scaled 1315    
\font\tenmsa=msam10          \font\ninemsa=msam9
\font\sevenmsa=msam7         \font\sixmsa=msam6
\font\fivemsa=msam5
%
%
\newfam\msafam                  \textfont\msafam=\tenmsa
\scriptfont\msafam=\sevenmsa    \scriptscriptfont\msafam=\fivemsa
\edef\hexa{\hexnumber\msafam}        
\def\msa{\fam\msafam\tenmsa}         
%
%
\font\thirtnmsb=msbm10 scaled 1315   
\font\tenmsb=msbm10      \font\ninemsb=msbm9
\font\sevenmsb=msbm7     \font\sixmsb=msbm6
\font\fivemsb=msbm5
%
\newfam\msbfam                   \textfont\msbfam=\tenmsb       
\scriptfont\msbfam=\sevenmsb     \scriptscriptfont\msbfam=\fivemsb
\edef\hexb{\hexnumber\msbfam}    
\def\msb{\fam\msbfam\tenmsb}     
%
%
\font\thirtneufm=eufm10 scaled 1315   
\font\teneufm=eufm10                 \font\nineeufm=eufm9
\font\seveneufm=eufm7                \font\sixeufm=eufm6
\font\fiveeufm=eufm5
%
\newfam\eufmfam                    \textfont\eufmfam=\teneufm
\scriptfont\eufmfam=\seveneufm     \scriptscriptfont\eufmfam=\fiveeufm
\edef\hexf{\hexnumber\eufmfam}      
\def\frak{\fam\eufmfam\teneufm}     
%
%
%
\font\thirtnrm=cmr10 scaled 1315    
\font\ninerm=cmr9                   \font\sixrm=cmr6   
%
\font\thirtni=cmmi10 scaled 1315    
\font\ninei=cmmi9                   \font\sixi=cmmi6  
%
\font\thirtnsy=cmsy10 scaled 1315   
\font\ninesy=cmsy9                  \font\sixsy=cmsy6  
%
\font\thirtnbf=cmbx10 scaled 1315   
\font\ninebf=cmbx9                  \font\sixbf=cmbx6  
%
%
\font\thirtnex=cmex10 scaled 1315   
\font\nineex=cmex9                  
%
%
\font\thirtnit=cmti10 scaled 1315  
\font\nineit=cmti9                  
%
\font\thirtnsl=cmsl10 scaled 1315  
\font\ninesl=cmsl9                  
%
\font\thirtntt=cmtt10 scaled 1315  
\font\ninett=cmtt9                  
%
%
%
%
\def\small{%
%
%
\textfont0=\ninerm \scriptfont0=\sixrm \scriptscriptfont0=\fiverm
\def\rm{\fam0\ninerm}
%
%
\textfont1=\ninei \scriptfont1=\sixi \scriptscriptfont1=\fivei
%
%
\textfont2=\ninesy \scriptfont2=\sixsy \scriptscriptfont2=\fivesy
%
%
\textfont3=\nineex \scriptfont3=\nineex \scriptscriptfont3=\nineex
%
%
\textfont\bffam=\ninebf \scriptfont\bffam=\sixbf
\scriptscriptfont\bffam=\fivebf \def\bf{\fam\bffam\ninebf}%
%
%
\textfont\itfam=\nineit \def\it{\fam\itfam\nineit}%
\textfont\slfam=\ninesl \def\sl{\fam\slfam\ninesl}%
\textfont\ttfam=\ninett \def\tt{\fam\ttfam\ninett}%
%
%
%
\textfont\msafam=\ninemsa \scriptfont\msafam=\sixmsa
\scriptscriptfont\msafam=\fivemsa \def\msa{\fam\msafam\ninemsa}%
%
%
\textfont\msbfam=\ninemsb \scriptfont\msbfam=\sixmsb
\scriptscriptfont\msbfam=\fivemsb \def\msb{\fam\msbfam\ninemsb}%
%
%
\textfont\eufmfam=\nineeufm  \scriptfont\eufmfam=\sixeufm
\scriptscriptfont\eufmfam=\fiveeufm \def\frak{\fam\eufmfam\nineeufm}%
%
%
%
\normalbaselineskip=11pt%
\setbox\strutbox=\hbox{\vrule height8pt depth3pt width0pt}%
%
%
\normalbaselines\rm
%
%
\stdskip=4pt plus2pt minus2pt    
\medskipamount=\stdskip          
\parskip=\stdskip                
\abovedisplayskip=\stdskip       
\belowdisplayskip=\stdskip       
\def\ppar{\par\goodbreak\vskip 6pt plus 3pt minus 3pt}%
%
%
\def\section##1{\global\advance\sectionnumber by 1
\vskip-\lastskip\penalty-800\vskip 20pt plus10pt minus5pt 
\egroup{\bf\number\sectionnumber\quad##1}\bgroup\small         
\vskip 6pt plus3pt minus3pt
\nobreak\resultnumber=1}
}    
%
\def\beginsmall{\bgroup\small}
\let\endsmall\egroup
%
%
%
%
\def\large{%
\textfont0=\thirtnrm \scriptfont0=\ninerm \scriptscriptfont0=\sevenrm
\def\rm{\fam0\thirtnrm}%
\textfont1=\thirtni \scriptfont1=\ninei \scriptscriptfont1=\seveni
\textfont2=\thirtnsy \scriptfont2=\ninesy \scriptscriptfont2=\sevensy
\textfont3=\thirtnex \scriptfont3=\thirtnex \scriptscriptfont3=\thirtnex
\textfont\bffam=\thirtnbf \scriptfont\bffam=\ninebf
\scriptscriptfont\bffam=\sevenbf \def\bf{\fam\bffam\thirtnbf}%
\textfont\itfam=\thirtnit \def\it{\fam\itfam\thirtnit}%
\textfont\slfam=\thirtnsl \def\sl{\fam\slfam\thirtnsl}%
\textfont\ttfam=\thirtntt \def\tt{\fam\ttfam\thirtntt}%
\textfont\msafam=\thirtnmsa \scriptfont\msafam=\ninemsa
\scriptscriptfont\msafam=\sevenmsa \def\msa{\fam\msafam\thirtnmsa}%
\textfont\msbfam=\thirtnmsb \scriptfont\msbfam=\ninemsb
\scriptscriptfont\msbfam=\sevenmsb \def\msb{\fam\msbfam\thirtnmsb}%
\textfont\eufmfam=\thirtneufm  \scriptfont\eufmfam=\nineeufm
\scriptscriptfont\eufmfam=\seveneufm \def\frak{\fam\eufmfam\teneufm}%
\normalbaselineskip=16pt%
\setbox\strutbox=\hbox{\vrule height11.5pt depth4.5pt width0pt}%
\normalbaselines\rm}%
%
%
\def\Bbb#1{{\msb#1}}

%

\def\re{\Bbb R}
%
\mathchardef\plussquare="0\hexa01
\mathchardef\nge="3\hexb0B
\mathchardef\maltesecross="0\hexa7A
\mathchardef\del="0\hexf01
%
%
%
%
\font\sc=cmcsc10
%
%
%
%
\def\sqr#1#2{{\vcenter{\vbox{\hrule  height.#2truept
	\hbox{\vrule width.#2truept height#1truept 
	\kern#1truept \vrule width.#2truept}
	\hrule height.#2truept}}}}
\def\sq{\sqr55}    
%
%
%
%
\newcount\sectionnumber            
\newcount\resultnumber             
\sectionnumber=0\resultnumber=1    
%
%
%
\def\section#1{\global\advance\sectionnumber by 1
\xdef\nextkey{\number\sectionnumber}
\vskip-\lastskip\penalty-800\vskip 20pt plus10pt minus5pt 
{\large\bf\number\sectionnumber\quad#1}         
\vskip 8pt plus4pt minus4pt
\nobreak\resultnumber=1}      
%
%
%
%
%
\def\sh#1{\vskip-\lastskip\ppar{\bf #1}\par\nobreak\medskip}         
%
%
%
%

%
\def\proc#1{\xdef\nextkey{\number\sectionnumber.\number\resultnumber}%
\vskip-\lastskip\ppar\bf%
\noindent#1\ \number\sectionnumber.\number\resultnumber
\stdspace\sl\global\advance\resultnumber by 1\ignorespaces}
\def\endproc{\rm\ppar} 
%
%
\def\prf{\vskip-\lastskip\ppar\noindent{\bf Proof}%
\stdspace\rm}                            
\def\endprf{\unskip\stdspace\hbox{}
\hfill$\sq$\par\medskip}                 
%
%
%
%
%
%
%
%
\def\proclaim#1{\vskip-\lastskip\ppar\bf%
\noindent#1\stdspace\sl\ignorespaces} 

%
%
%
%
\def\rk#1{\vskip-\lastskip\ppar{\bf #1}\stdspace\ignorespaces}                
\def\endrk{\par\medskip}
%
%
%
%
%
%
\def\label{\xdef\nextkey{\number\sectionnumber.\number\resultnumber}%
\number\sectionnumber.\number\resultnumber
\global\advance\resultnumber by 1}
%
%
%
%
%
%
%
%
%
%
%
%
%
%
%
%
\newcount\refnumber              
\refnumber=1                     
\long\def\reflist#1\endreflist{%
\long\def\thereflist{#1}{\def\refkey##1##2\par{\xdef##1{\number\refnumber}%
\global\advance\refnumber by 1}%
\def\key##1##2\par{\expandafter\xdef%
\csname##1\endcsname{\number\refnumber}%
\global\advance\refnumber by 1}#1\par}}
\long\def\references{%
\penalty-800\vskip-\lastskip\vskip 15pt plus10pt minus5pt 
{\large\bf References}\ppar 
{\leftskip=25pt\frenchspacing    
\small\parskip=3pt plus2pt       
\def\refkey##1##2\par{\noindent  
\llap{[##1]\stdspace}\ignorespaces##2\par}         
\def\key##1##2\par{\noindent  
\llap{[\ref{##1}]\stdspace}\ignorespaces##2\par}  
\def\,{\thinspace}\thereflist\par}}
%
%
%
\newcount\footnotenumber         
\footnotenumber=1                
\def\fnote#1{\xdef\nextkey{\number\footnotenumber}%
{\small\ifnum\footnotenumber>9\parindent=14pt%
\else\parindent=10pt\fi\footnote{$^{\number\footnotenumber}$}%
{\hglue-5pt#1}\global\advance\footnotenumber by 1}}
%
%
%
%
%
%
%
\newcount\figurenumber          
\figurenumber=1                 
\def\caption#1{\xdef\nextkey{\number\figurenumber}%
\cl{\small Figure \number\figurenumber: #1}%
\global\advance\figurenumber by 1}
\def\figurelabel{\xdef\nextkey{\number\figurenumber}%
\cl{\small Figure \number\figurenumber}%
\global\advance\figurenumber by 1}
\long\def\figure#1\endfigure{{\xdef\nextkey{\number\figurenumber}%
\let\captiontext\relax\def\caption##1{\xdef\captiontext{##1}}%
\midinsert\cl{\ignorespaces#1\unskip\unskip\unskip\unskip}\vglue6pt\cl{\small 
Figure \number\figurenumber\ifx\captiontext\relax\else: \captiontext
\fi}\endinsert\global\advance\figurenumber by 1}}
%
%
%
%
%
%
%
\def\nextkey{??}   
%
\def\key#1{\expandafter\xdef\csname #1\endcsname{\nextkey}}
\def\ref#1{\expandafter\ifx\csname #1\endcsname\relax
\immediate\write16{Reference {#1} undefined}??\else
\csname #1\endcsname\fi}
%
%
%
%
%
%
%
\newread\gtinfile
\newwrite\gtreffile
\def\useforwardrefs{
\openin\gtinfile\jobname.ref
\ifeof\gtinfile
\closein\gtinfile
\immediate\write16{No file \jobname.ref}
\else
\closein\gtinfile
\input \jobname.ref
\fi
\immediate\openout\gtreffile \jobname.ref
%
%
\def\key##1{{\def\\{\noexpand}%
\expandafter\xdef\csname ##1\endcsname{\nextkey}%
\immediate\write\gtreffile{\\\expandafter\\\def\\\csname ##1\\\endcsname%
{\nextkey}}}}
%
%
\long\def\reflist##1\endreflist{%
\long\def\thereflist{##1}{\def\refkey####1####2\par{\xdef####1{%
\number\refnumber}{\def\\{\noexpand}\immediate\write\gtreffile
{\\\def\\####1{\number\refnumber}}}\global\advance\refnumber by 1}%
\def\key####1####2\par{\expandafter\xdef%
\csname####1\endcsname{\number\refnumber}%
{\def\\{\noexpand}\immediate\write\gtreffile
{\\\expandafter\\\def\\\csname ####1\\\endcsname{\number\refnumber}}}
\global\advance\refnumber by 1}##1\par}}
\long\def\biblio##1\endbiblio{\reflist##1\endreflist\references}%
%
%
\def\numkey##1{{\def\\{\noexpand}%
\xdef##1{\number\sectionnumber.\number\resultnumber}
\immediate\write\gtreffile{\\\def\\##1%
{\number\sectionnumber.\number\resultnumber}}}}
\def\seckey##1{{\def\\{\noexpand}\xdef##1{\number\sectionnumber}
\immediate\write\gtreffile{\\\def\\##1{\number\sectionnumber}}}}
\def\figkey##1{\xdef##1{\number\figurenumber}%
{\def\\{\noexpand}\immediate\write\gtreffile%
{\\\def\\##1{\number\figurenumber}}}
\number\figurenumber\global\advance\figurenumber by 1}
}   
%
%
%
%
\def\figkey#1{\xdef#1{\number\figurenumber}%
\number\figurenumber\global\advance\figurenumber by 1}
\def\fig#1#2\endfig{%
\midinsert\cl{#2}\vglue6pt\cl{\small Figure #1}\endinsert}
\def\newfig{\number\figurenumber\global\advance\figurenumber by 1}
\def\numkey#1{\xdef#1{\number\sectionnumber.\number\resultnumber}}
\def\seckey#1{\xdef#1{\number\sectionnumber}}
%
%
%
%
%
%
%
%
%
\def\verb{\catcode`\"=\active}       
\def\brev{\catcode`\"=12}            
\brev                                
\verb                                
{\obeyspaces\gdef {\ }}              
{\catcode`\`=\active\gdef`{\relax\lq}}
\def"{%
\begingroup\baselineskip=12pt\def\par{\leavevmode\endgraf}%
\tt\obeylines\obeyspaces\parskip=0pt\parindent=0pt%
\catcode`\$=12\catcode`\&=12\catcode`\^=12\catcode`\#=12%
\catcode`\_=12\catcode`\~=12%
\catcode`\{=12\catcode`\}=12\catcode`\%=12\catcode`\\=12%
\catcode`\`=\active\let"\endgroup}
\brev      
%
%
%
%
%
%
\def\items{\par\leftskip = 25pt}           
\def\enditems{\par\leftskip = 0pt}         
\def\item#1{\par\leavevmode\llap{#1\stdspace}%
\ignorespaces}                             
%
%

%
%
\def\co{\colon\thinspace}    
\def\np{\vfil\eject}         
\def\nl{\hfil\break}         
\def\cl{\centerline}         
%
%
%

%
%
%
%
%
\def\title#1{\def\thetitle{#1}}

\def\author#1{\edef\previousauthors{\theauthors}
 \ifx\theauthors\relax\def\theauthors{#1}\else
 \def\theauthors{\previousauthors\par#1}\fi}

\let\authors\author        
\def\address#1{\edef\previousaddresses{\theaddress}
 \ifx\theaddress\relax\def\theaddress{#1}\else
 \def\theaddress{\previousaddresses\par\vskip 2pt\par#1}\fi}
\def\secondaddress#1{\edef\previousaddresses{\theaddress}
 \ifx\theaddress\relax\def\theaddress{#1}\else
 \def\theaddress{\previousaddresses\par{\rm and}\par#1}\fi}   

\def\email#1{\edef\previousemails{\theemail}
 \ifx\theemail\relax\def\theemail{#1}\else
 \def\theemail{\previousemails\hskip 0.75em\relax#1}\fi}
\def\secondemail#1{\edef\previousemails{\theemail}
 \ifx\theemail\relax\def\theemail{#1}\else
 \def\theemail{\previousemails\hskip 0.75em{\rm and}\hskip 0.75em
 \relax#1}\fi}
\def\url#1{\edef\previousurls{\theurl}
 \ifx\theurl\relax\def\theurl{#1}\else
 \def\theurl{\previousurls\hskip 0.75em\relax#1}\fi}
\def\secondurl#1{\edef\previousurls{\theurl}
 \ifx\theurl\relax\def\theurl{#1}\else
 \def\theurl{\previousurls\hskip 0.75em{\rm and}\hskip 0.75em
 \relax#1}\fi}
\long\def\abstract#1\endabstract{\long\def\theabstract{#1}}
\def\primaryclass#1{\def\theprimaryclass{#1}}
\def\secondaryclass#1{\def\thesecondaryclass{#1}}
\def\keywords#1{\def\thekeywords{#1}}
%
%
\let\\\par\let\thetitle\relax
\let\theauthors\relax
\let\theaddress\relax\let\theshortaddress\relax
\let\theemail\relax\let\theurl\relax
\let\theabstract\relax\let\theprimaryclass\relax
\let\thesecondaryclass\relax\let\thekeywords\relax
%
%
%
%
\long\def\maketitlepage{    

\vglue 0.2truein   

%
{\parskip=0pt\leftskip 0pt plus 1fil\def\\{\par\smallskip}{\large
\bf\thetitle}\par\medskip}   

\vglue 0.15truein 

%
{\parskip=0pt\leftskip 0pt plus 1fil\def\\{\par}{\sc\theauthors}
\par\medskip}%
 
\vglue 0.1truein 

%
{\small\parskip=0pt
{\leftskip 0pt plus 1fil\def\\{\par}{\sl\theaddress}\par}
\ifx\theemail\relax\else  
\vglue 5pt \def\\{\stdspace{\rm and}\stdspace} 
\cl{Email:\stdspace\tt\theemail}\fi
\ifx\theurl\relax\else    
\vglue 5pt \def\\{\stdspace{\rm and}\stdspace} 
\cl{URL:\stdspace\tt\theurl}\fi\par}

\vglue 7pt 

{\bf Abstract}

\vglue 5pt

\theabstract

\vglue 7pt 

{\bf AMS Classification numbers}\quad Primary:\quad \theprimaryclass\par

Secondary:\quad \thesecondaryclass

\vglue 5pt 

{\bf Keywords:}\quad \thekeywords

\np  

}    
%
%
\long\def\makeshorttitle{    


%
{\parskip=0pt\leftskip 0pt plus 1fil\def\\{\par\smallskip}{\large
\bf\thetitle}\par\medskip}   

\vglue 0.05truein 

%
{\parskip=0pt\leftskip 0pt plus 1fil\def\\{\par}{\sc\theauthors}
\par\medskip}%
 
\vglue 0.03truein 

%
{\small\parskip=0pt
{\leftskip 0pt plus 1fil\def\\{\par}{\sl\ifx\theshortaddress\relax
\theaddress\else\theshortaddress\fi}\par}
\ifx\theemail\relax\else  
\vglue 5pt \def\\{\stdspace{\rm and}\stdspace} 
\cl{Email:\stdspace\tt\theemail}\fi
\ifx\theurl\relax\else    
\vglue 5pt \def\\{\stdspace{\rm and}\stdspace} 
\cl{URL:\stdspace\tt\theurl}\fi\par}

\vglue 10pt 


{\small\leftskip 25pt\rightskip 25pt{\bf Abstract}\stdspace\theabstract

{\bf AMS Classification}\stdspace\theprimaryclass
\ifx\thesecondaryclass\relax\else; \thesecondaryclass\fi\par
{\bf Keywords}\stdspace \thekeywords\par}
\vglue 7pt
}    
\let\maketitle\makeshorttitle      
\input epsf
%
%
\def\relabelbox{%
  \hbox\bgroup%
}%
\def\endrelabelbox{%
\egroup%
}%
\def\relabel #1#2 {%
  \special{ps:/a {} def}%
  \smash{\rlap{#2}}%
}%
\def\adjustrelabel <#1,#2> #3#4 {%
  \special{ps:/a {} def}%
  \smash{\rlap{\kern #1 \raise #2\hbox{#4}}}%
}%
\def\extralabel <#1,#2> #3 {\smash{\rlap{\kern #1 \raise #2\hbox{#3}}}}%

\chardef\newinsCatAt\the\catcode `\@
\catcode `\@=11
%
%
%
\newskip\insertskipamount\newskip\inserthardskipamount
\insertskipamount 12pt plus2pt  
\inserthardskipamount 4pt       
\def\insertskip{\vskip\insertskipamount}
%
%
\newskip\LastSkip
\def\SaveLastSkip{\LastSkip\lastskip}
\def\RestoreLastSkip{\nobreak\vskip-\LastSkip\vskip\LastSkip}
%
%
\newcount\SplitTest
\def\SetSplitTest{\SplitTest\insertpenalties
  \insert\topins{\floatingpenalty1}%
  \advance\SplitTest-\insertpenalties}
%
%
\def\midinsert{\par
 \SaveLastSkip\penalty-150\SetSplitTest\RestoreLastSkip
 \ifnum\SplitTest=-1
  \@midfalse\p@gefalse\else\@midtrue\fi\@ins}
\def\@ins{\par\begingroup\setbox\z@\vbox\bgroup%
  \vglue\inserthardskipamount}
\def\endinsert{\egroup 
  \if@mid \dimen@\ht\z@ \advance\dimen@\dp\z@
    \advance\dimen@\insertskipamount
    \advance\dimen@\pagetotal\advance\dimen@-\pageshrink
    \ifdim\dimen@>\pagegoal\@midfalse\p@gefalse\fi\fi
  \if@mid%
    \ifdim\lastskip<\insertskipamount\removelastskip\insertskip\fi
    \nointerlineskip\box\z@\penalty-200\insertskip
  \else%
    \SaveLastSkip
    \insert\topins{\penalty100 
    \splittopskip\z@skip
    \splitmaxdepth\maxdimen \floatingpenalty\z@
    \ifp@ge \dimen@\dp\z@
    \vbox to\vsize{\unvbox\z@\kern-\dimen@}
    \else \box\z@\nobreak\insertskip\fi}
    \RestoreLastSkip
   \fi\endgroup}
%
\catcode `\@=\newinsCatAt
 
%
\hoffset 0.5truein     
\voffset 1truein       
%
%
%
\def\tw{{\tt \char'176}}  

\def\ep{\varepsilon}

\font\spec=cmtex10 scaled 1095 
\def\d{\hbox{\spec \char'017\kern 0.05em}}

\def\re{\Bbb R}
\def\Z{\Bbb Z}

\def\ss{\scriptstyle}

\reflist

\key{Carter}
{\bf J Scott Carter}, {\bf Daniel Jelsovsky}, {\bf Laurel Langford},
{\bf Seichi Kamada}, {\bf Masahico Saito}, {\it Quandle Cohomology and
State-sum Invariants of Knotted Curves and Surfaces}, 
{\tt arxiv:math.GT/9906115}

\key{Racks}
{\bf Roger Fenn}, {\bf Colin Rourke}, {\it Racks and links in
codimension two}, Journal of Knot theory and its Ramifications, 1
(1992) 343--406, available from:\nl {\tt
http://www.maths.warwick.ac.uk/\tw cpr/ftp/racks.ps}

\key{Trunks}
{\bf Roger Fenn}, {\bf Colin Rourke}, {\bf Brian Sanderson},
{\it Trunks and classifying spac\-es}, Applied categorical structures,
3 (1995) 321--356, available from:\nl {\tt
http://www.maths.warwick.ac.uk/\tw cpr/ftp/trunks.ps}

\key{James}
{\bf Roger Fenn}, {\bf Colin Rourke}, {\bf Brian Sanderson}, {\it
James bundles and applications}, Warwick preprint (1996), available
from:\nl {\tt http://www.maths.warwick.ac.uk/\tw cpr/ftp/james.ps}

\key{Hempel}
{\bf John Hempel}, {\it 3--manifolds}, Annals of Math. Study no. 86, Princeton
University Press (1976)

\key{Mono}
{\bf Roger Fenn}, {\bf Colin Rourke}, {\bf Brian Sanderson}, {\it
The classification of links}, Monograph in preparation

\key{Comp}
{\bf Colin Rourke}, {\bf Brian Sanderson}, {\it The compression
theorem I}, Geometry and Topology, 5 (2001) 399--429

\key{Maple}
{\bf Brian Sanderson}, {\it Maple worksheet}, 
\nl {\tt
http://www.maths.warwick.ac.uk/\tw bjs/trefoil.ms}

\endreflist

\title{A new classification of links and\\some calculations using it}

\authors{Colin Rourke\\Brian Sanderson}

\address{Mathematics Institute, University of Warwick\\Coventry, CV4 7AL, UK}

\email{cpr@maths.warwick.ac.uk\\bjs@maths.warwick.ac.uk}
\url{http://www.maths.warwick.ac.uk/\tw cpr/\qua {\rm and}\qua \tw bjs/}

\abstract
A new classification theorem for links by the authors and Roger Fenn
leads to computable link invariants. As an illustration we distinguish
the left and right trefoils and recover the result of Carter et al
that the 2-twist-spun trefoil is not isotopic to its orientation
reverse. We sketch the proof the classification theorem. Full details
will appear elsewhere

\endabstract

\primaryclass{57Q45}\secondaryclass{57M25, 57M27}


\makeshorttitle

\section{Introduction}

Throughout the paper, links and link diagrams will be assumed to be
framed, in other words to have trivialised normal bundle.

In [\Trunks, \James] Fenn, Rourke and Sanderson define a space $BL$
associated to a link $L$ in codimension 2 which is the {\it rack
space\/} of the fundamental rack of the link.  The link determines a
{\it canonical class\/} in the homotopy of $BL$.  For example, if the
link is in $S^3$ then the canonical class $c(L)$ is in $\pi_2(BL)$.

Furthermore in [\Carter] Carter, Jelsovsky, Kamada, Langford and Saito
construct ``state-sum'' invariants of knots using a notion of quandle
cohomology.  As an application they prove that the twice twist-spun
trefoil is not isotopic to its orientation reverse.  Quandle
cohomology groups are closely related to the cohomology groups of the
rack space, and the state-sum invariant has a natural interpretation
in terms of the canonical class.  This suggests the importance of this
canonical class.

In this paper we announce, with outline proof, a new result: the
fundamental rack together with the canonical class classifies
classical links.

\proclaim{Classification Theorem}{\rm (Fenn--Rourke--Sanderson)}\qua
Suppose that $L,M$ are two links in $S^3$ and suppose that there is an
isomorphism of fundamental racks $\phi\co \Gamma(L)\to \Gamma(M)$ such
that $\phi_*(c(L))=c(M)$ then $L$ and $M$ are isotopic.\endproc

This result should be contrasted with existing classification results
using the fundamental rack.  The augmented fundamental rack classifies
{\it irreducible\/} links in a 3--manifold up to homeomorphism
[\Racks; 5.2] and the unaugmented fundamental rack does the same for
irreducible links in a homotopy 3--sphere (and in particular for
classical links) [\Racks; 5.3].  However because the fundamental rack
of the {\it inverse mirror\/} of a link (meaning the mirror image with
orientation change) is isomorphic to the original link, the
fundamental rack alone cannot classify links up to isotopy.  For
example the fundamental racks of the left and right-handed trefoil
knots are isomorphic.  The situation is compounded for reducible links
where the fundamental rack only determines the link up to inversion of
the {\it blocks} (maximal irreducible sublinks).

The canonical class encodes the orientations of all the blocks and
together with the fundamental rack enables a complete reconstruction
of the link.  It follows that all link invariants can theoretically be
recovered from the rack and canonical class.  Now the canonical class
lies in $\pi_2$ of the rack space which is $H_2$ of the universal
cover.  Partial information can be obtained from $H_2$ of a cover
intermediate between the rack space and the universal cover and by
using a representation in a finite rack, this homology group becomes
computable.  Thus we have a simple way to construct computable link
invariants which approximate in a systematic way to the full invariant
of rack and canonical class.

In this paper we will demonstrate the practicality of this by
constructing a computable invariant which distinguishes the left and
right-handed trefoils.  The calculation needed is of the third
homology group of the three-colour rack and for this we rely on a
Maple worksheet [\Maple].  The same calculation allows us to recover
the theorem of Carter et al [\Carter].

The Classification Theorem is joint work with Roger Fenn.  A full
treatment will be given in [\Mono].  Here we prove the result using
results from [\James].

Here is an outline of the paper.  In section 2 we give some background
material on racks and the rack space and prove the Classification
Theorem.  In section 3 we start our investigation of the use of the
canonical class by using it together with a short Maple calculation to
prove that the 2--twist-spun trefoil is not isotopic to its
orientation reverse preserving orientations.  In section 4 we extend
the notion of the rack space and use the same calculation to prove
that the left and right trefoils are different and we conclude in
section 5 with some remarks about the proofs.

\section{Racks and the Classification Theorem}

\sh{The rack space and the canonical class of a diagram}

Full details of all the material in this section can be found in
[\James].  For background material on racks see [\Racks].  Here we
shall summarise the results from [\James] which are needed for the
Classification Theorem and the later calculations.

Let $R$ be a rack.  The {\sl rack space} $BR$ of $R$ is a cubical set
with the set of $n$--cubes in bijection with $R^n$.  In this paper we
shall only need to consider cubes of dimension $\le3$.  The
interpretation of such cubes is given by the pictures in figure 1.

\figure
\relabelbox\epsfxsize 5truein
\epsfbox{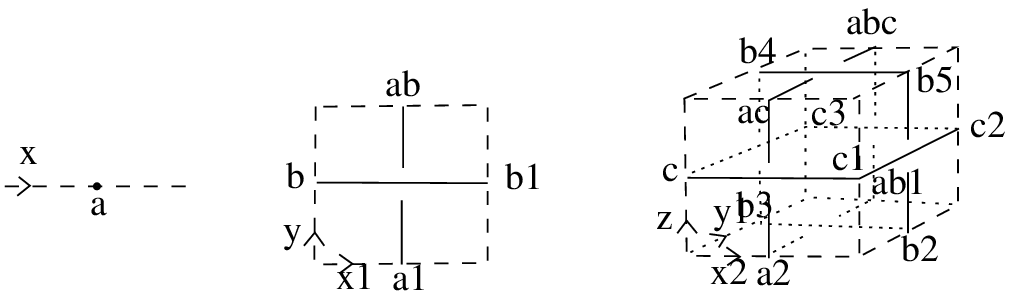}
\relabel {x}{\small $\ss x$}
\relabel {x1}{\small $\ss x$}
\relabel {x2}{\small $\ss x$}
\relabel {y}{\small $\ss y$}
\relabel {y1}{\small $\ss y$}
\relabel {z}{\small $\ss z$}
\relabel {a}{\small $a$}
\relabel {a1}{\small $a$}
\relabel {a2}{\small $a$}
\relabel {ab}{\small $a^b$}
\relabel {b}{\small $b$}
\relabel {b1}{\small $b$}
\relabel {b2}{\small $b$}
\relabel {b4}{\small $b^c$}
\relabel {b5}{\small $b^c$}
\relabel {c}{\small $c$}
\relabel {c1}{\small $c$}
\relabel {c2}{\small $c$}
\relabel {ac}{\small $a^c$}
\relabel {abc}{\small $a^{bc}=a^{cb^c}$}
\relabel {c3}{$\ss c$}
\relabel {ab1}{$\ss a^b$}
\relabel {b3}{$\ss b$}
\endrelabelbox
\endfigure

The figures make the boundary maps clear.  Thus the back face of the
3--cube drawn $(a,b,c)$ is glued to the 2--cube $(a^b,c)$ and the top
face is glued to $(a^c,b^c)$.

There are similar descriptions of boundary maps from $n$--cubes to
$(n-1)$--cubes and the general formula is:
$$\eqalign{\d_i^0(x_1,\ldots ,x_n)&=(x_1,\ldots
,x_{i-1},x_{i+1},\ldots ,x_n),\cr
\d_i^1(x_1,\ldots ,x_n)&=
(x_1^{x_i},\ldots ,x_{i-1}^{x_i},x_{i+1},\cdots
,x_n)
\hbox {\quad for\quad}1\leq i\leq n.}$$
The rack space determines a topological space (also denoted $BR$) by
gluing real cubes together using these boundary maps.

Now suppose that we are given a diagram in $\re^2$ representing a link
$L$ in $\re^3$.  A {\sl labelling} of the diagram in a rack $R$ means
a labelling of arcs by elements of $R$ so that at double points the
rule indicated in figure 2 holds.

\figure
\relabelbox\small
\epsfbox{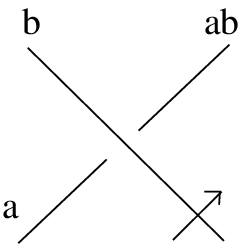}
\relabel {a}{$a$}
\relabel {ab}{$a^b$}
\relabel {b}{$b$}
\endrelabelbox
\endfigure

A labelling by $R$ is precisely the same as a homomorphism of the
fundamental rack of $L$ to $R$, see [\Racks; page 384].

Now $L$ is canonically framed by taking the first framing vector
parallel to the vertical and second vector normal to the image in
$\re^2$.  Using this second vector we can provide the diagram with a
bicollar with collar lines identified with $I$ (so that the
submanifold corresponds to ${1\over2}\in I$ and the vector points in
the positive direction along $I$).  We can assume that the collar
lines determine little squares at double points.  The squares are
canonically identified with $I^2$ using the vertical ordering of arcs.
Now suppose the diagram is labelled in $R$ then using the labelling
these copies of $I$ and $I^2$ can be identified with cubes in $BR$.
The diagram now determines a map of $\re^2$ to $BR$ by mapping the
collar lines to the 1--cubes of $BR$ with which they are identified
and the little squares at double points to the appropriate 2--cubes.
Outside the collar all is mapped to the basepoint (the unique 0--cube)
and hence we get a based map $S^2\to BR$ by mapping infinity (the
basepoint of $S^2$) to the basepoint of $BR$.

This map is almost canonical, depending only on the choice of collar,
and we call it the {\sl canonical map} determined by the diagram and
the corresponding class in $\pi_2(BR)$ the {\sl canonical class}.
Given an isotopy of diagrams (or even a bordism) labelled in $R$,
then the canonical map varies by a homotopy by applying a similar
construction to the bordism.

There is a converse construction.  Given a map $S^2\to BR$ we can use
transversality to make it transverse to the 2--skeleton of $BR$ and
this means that it pulls back a diagram in $\re^2$ (for some link in
$\re^3$) labelled in $R$.  Similarly a homotopy pulls back a bordism
of diagrams, ie, a diagram for a bordism emdedded in $\re^3\times I$.

These considerations are summarised in the following theorem:

\proc{Theorem}\key{RackClass}There is a bijection between $\pi_2(BR)$ 
and the set of bordism classes of diagrams in $\re^2$ labelled in
$R$.\endproc

See [\James] for the detailed proof of this result, for the
interpretation in terms of James bundles and for more general results
along the same lines.

\sh{Proof of the Classification Theorem}

By theorem \ref{RackClass} the data in the theorem determines a
bordism $W$ between $L$ and $M$ with homomorphism $\Phi$ from
$\Gamma(W)$ to $\Gamma(L)$ (ie a labelling in $\Gamma(L)$) which is
the identity on $L$ and on $M$ is the given isomorphism $\phi$.

We shall replace $W$ by an isotopy (ie a level-preserving embedding of
a link $\times I$) by induction on the number of blocks (ie maximal
irreducible sublinks).  Suppose that $L = L_1 \coprod L_2$ where $L_1$
is a block.  Then the fundamental rack of $L$ is the free product of
the fundamental racks of $L_1$ and $L_2$ and likewise the fundamental
group (see [\Racks; page 357] for the notion of the free product of
two racks).  Now the fundamental group is the associated group of the
fundamental rack (this is true for classical links by [\Racks; 3.3]
--- it is not true for general links) and it follows that the
fundamental group of $M$ is also a free product.  Then by Kneser's
conjecture [\Hempel; pages 66--68] $M$ is a split link and furthermore
the proof shows that we can assume that the splitting induces the same
free product decomposition of the fundamental group.  We shall now
relate the splittings of $L$ as $L_1 \coprod L_2$ and of $M$, just
constructed.

Think of the components of $W$ as coloured red if they are labelled by
elements of $\Gamma(L)$ in one of the orbits of $L$ corresponding to
$L_1$ (recall that orbits in a fundamental rack of a link are in
bijection with components of the link) and coloured blue of they are
likewise labelled in orbits corresponding to $L_2$.  We obtain a
decomposition of $W$ as $W_1\cup W_2$ where $W_1$ is red and $W_2$ is
blue and an induced decomposition $M=M_1\cup M_2$ of $M$.  We claim that
this is the same as the splitting of $M$ constructed above.

To see this consider the augmentation $\d$ which takes a element of
the fundamental rack to the element of the fundamental group which is
the ``frying pan'': down the rack element (represented as a path from
the framing curve to the basepoint) round the meridian and back up the
rack element.  It follows that the augmentation takes an element in an
orbit (ie component) to a fundamental group element which links the
component once (and others zero times).  But an element of
$\Gamma(M_1)$ is mapped to an element in $\Gamma(L_1)$ and hence is
augmented to an element in the free factor of the fundamental group
corresponding to $L_1$, which was how we obtained the splitting of $M$
using Kneser.  Thus $M$ splits as $M_1\coprod M_2$ with $M_1$ coloured
red and $M_2$ coloured blue.  It follows that $\Gamma(M)$ is the free
product $\Gamma(M_1)*\Gamma(M_2)$.

Now we claim that the isomorphism $\phi$ splits as an isomororphism
respecting the two free product decompositions.  To see this we
observe that, given an element of a rack, its augmentation is
non-trivial in the fundamental group (it links the appropriate arc
once) and hence if we consider an element of $\Gamma(M_1)$ then its
augmentation is non-trivial in $\pi_1(M_1)$ and maps to $\pi_1(L_1)$
since the splittings of the groups correspond.  It follows that the
image under $\phi$ is in $\Gamma(L_1)$.  Similarly $\phi$ takes
$\Gamma(M_2)$ to $\Gamma(L_2)$.

Now we perform a trick.  Think of the splittings of $L$ and $M$ as
into left and right halves (red to the left, blue to the right).  $W$
is not split, but we pull it apart by translating the red pieces to
the left past the blue pieces and form in this way a new split
bordism $W'=W_1\coprod W_2$.  We replace the labels on $W_1$ by the
old labels on $W$ composed with the homomorphism
$\Gamma(L)=\Gamma(L_1)*\Gamma(L_2)\to \Gamma(L_1)\coprod\Gamma(L_2)$
given by ignoring the $L_2$ contributions to elements corresponding to
$L_1$ and vice versa (recall the description of the free product on
[\Racks; page 357]).  This gives a new labelling but observe that the
labelling on $M$ is undisturbed by this trick.

We now have new bordisms $W_1$ between $L_1$ and $M_1$ and $W_2$
between $M_1$ and $M_2$.  We now use the proof of [\James; 5.10].
There is no obstruction to constructing a map of $S^3\times I$ to
$S^3$ which is transverse to $L_1$ and such that the preimage of of
$L_1$ is $W_1$ (essentially because $S^3-L_1$ is a $K(\pi,1))$.  This
map is the identity on the $L_1$ end and on the $L_2$ and deforms to a
homeomorphism as in the proof of [\Racks; 5.2] (essentially this is
Waldhausen's theorem).  Thus we have a homotopy between the identity
and a homeomorphism of $L_1$ to $L_2$ which can be replaced by an
isotopy since homotopic homeomorphisms of spheres are isotopic.

Finally the other half of $W'$, namely $W_2$ can also be replaced by
an isotopy by induction and the theorem is proved.

\section{The 2--twist spun trefoil}

We start our calculations using the canonical class by considering
2--knots (ie knots of $S^2$ in $\re^4$).  This is because
distinguishing the left and right trefoils introduces extra technical
details.

\sh{Framings and diagrams for 2--knots in 4--space}

For classical knots, the oriented and framed theories differ, for
example the writhe of a knot is an invariant of framed knots but not
of oriented knots.  For 2--knots the distinction disappears:

Let $K\subset\re^4$ be an oriented, possibly knotted, 2--sphere in
$\re^4$.  The normal bundle of $K$ has a section since $K$ bounds a
3--manifold (a generalised ``Seifert surface'') and furthermore the
orientations now give a canonical second section, in other words a
framing of $K$.  Thus we can confuse sections of the normal bundle and
framings.  Now two sections of $K$ differ by a map $K\to S^1$ (since
the normal bundle is trivial) and since $\pi_2(S^1)=0$ any two
sections are isotopic.  Thus the orientations determine a canonical
framing of $K$ in $\re^4$.

Now consider the canonical projection $\re^4\to\re^3$ and think of
$\re^3$ as horizontal and the remaining $\re^1$ as vertical.  Suppose
that we are given a 2--knot $K$ such that the projection of $K$ on
$\re^3$ has multiple set comprising transverse double lines and
triple points.  The image of the projection together with the
information of the vertical order of sheets at the multiple set is
called a {\sl knot diagram} for $K$.  The diagram determines $K$ up to
isotopy and also determines a natural framing of $K$ corresponding to
the section which is vertically up.  Conversely, by the Compression
Theorem [\Comp] a 2--knot $K$ with a section of its normal bundle can
be isotoped so that the section is vertically up and hence $K$
projects to an immersion in $\re^3$.  By making this immersion
self-transverse we now have a knot diagram for $K$.  Furthermore, by
the 1--parameter version of the Compression Theorem, isotopic framed
embeddings determine isotopic knot diagrams.  Combining the two facts
just outlined, namely that orientation determines framing and that
framing determines diagram, we have:

\proc{Lemma}\key{Bijection}There is a natural bijection between isotopy 
classes of oriented embeddings of a 2--sphere in $\re^4$ and isotopy
classes of knot diagrams in $\re^3$.\endproc

By a similar construction to that given in section 2 for diagrams in
the plane, a diagram of a 2--knot whose fundamental rack is labelled
in a rack $R$ determines a canonical class in $\pi_3(BR)$.  Using
lemma \ref{Bijection} we have:

\proc{Lemma}\key{Canon}The
canonical class in $\pi_3(BR)$ of a labelled 2--knot diagram is
determined by\items
\item{\rm(a)}the isotopy class of the oriented knot
\item{\rm(b)}the labelling, or equivalently, the homomorphism of the 
fundamental rack to $R$.\enditems\endproc

There is a natural extension of theorem \RackClass\ to 2--knots:

\proc{Theorem}\key{RackClassBis}There is a bijection between $\pi_3(BR)$ 
and the set of bordism classes of diagrams in $\re^3$ of 2--knots
labelled in $R$.\endproc

We use a calculation involving the canonical class to prove:

\proc{Theorem}The 2--twist spun trefoil is not (oriented)
isotopic to its orientation reverse.\endproc

\prf The proof occupies the remainder of this section.  Let $K$ be 
the 2--twist spun trefoil.  We shall exhibit an explicit diagram for
$K$ which is labelled by the ``three-colour rack'' $T:=\{0,1,2\}$ with
$a^b:= 2b-a$ mod 3, ie, $a^b=c$ iff $a,b,c$ are all the same or all
different.

\figure\epsfxsize\hsize
\epsfbox{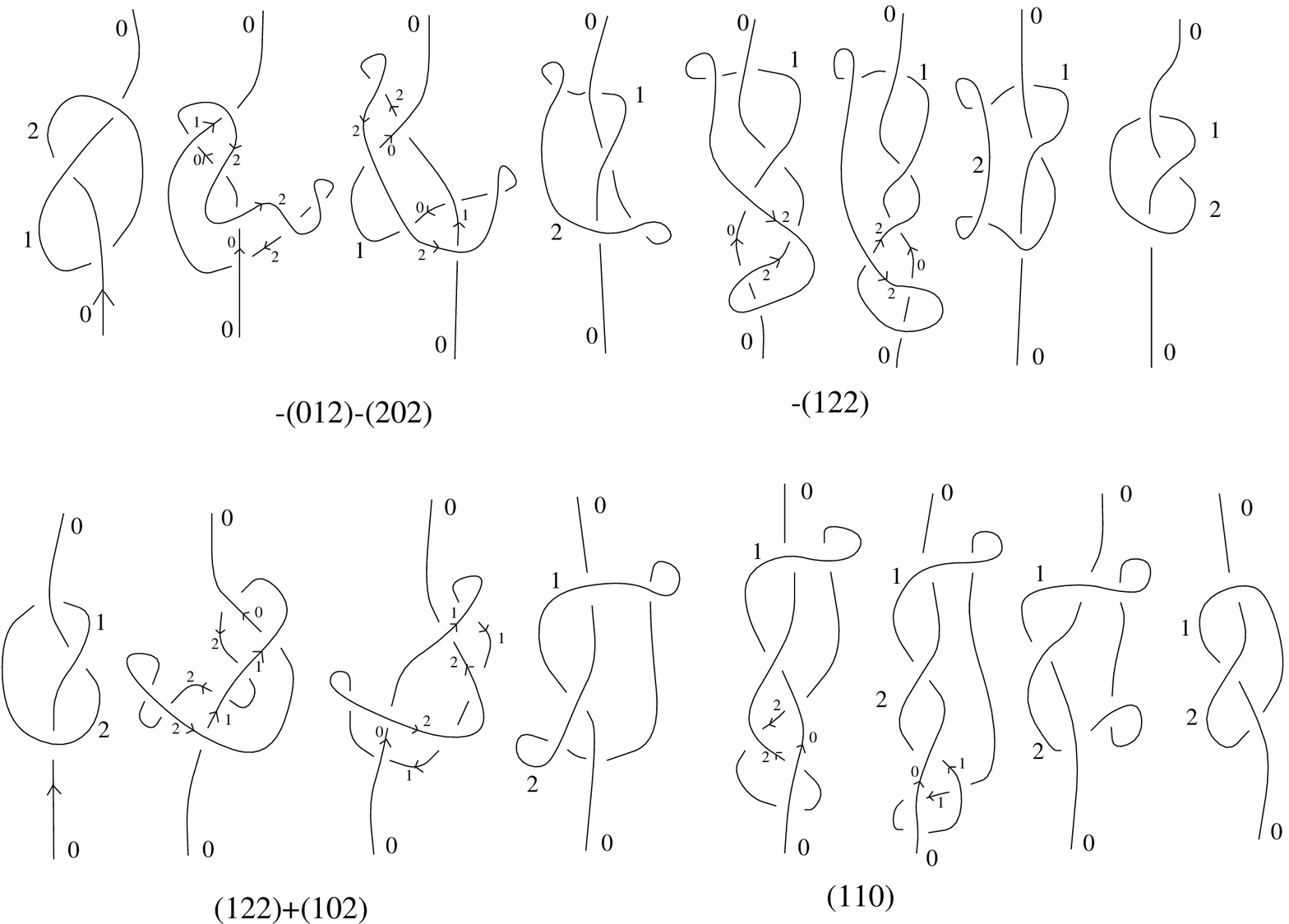}
\endfigure

In figure 3 we give a series of slices of the 2--twist spun trefoil.
The pictures are to be understood as follows.  We are thinking of
$\re^4$ as an ``open-book decomposition'' $\re^3_+\times S^1$ with
$x\times S^1$ identified with $x$ for each $x\in\re^2=\d\re^3_+$, ie,
$\re^2$ is the ``spine'' of the decomposition.  Similarly we are
thinking of $\re^3$ as an open-book decomposition with spine $\re^1$,
ie, $\re^3$ is $\re^2_+ \times S^1$ with a similar identification for
each $x\in\re^1=\d\re^2_+$.  The projection $\re^4\to\re^3$ is then
given by the standard projection $\re^3_+\to\re^2_+$ crossed with the
identity on $S^1$.  The figures are drawn as a sequence of diagrams in
copies of $\re^2_+$ corresponding to embeddings in the corresponding
copies of $\re^3_+$ and describe a diagram in $\re^3$ for an embedding
in $\re^4$.  The basepoint in $\re^4$ is chosen to lie on the spine
$\re^2$ of the decomposition above the spine $\re^1$ on the
decomposition of $\re^3$.  Thus the basepoint lies above each of the
diagrams drawn.  The framing (not shown) is given by the left-hand
rule in other words the framing vector is obtained from the
orientation vector on arcs by turning to the left.  Figure 3 shows one
twist of the spun trefoil.  To get the 2--twist spun trefoil, we
repeat the sequence.  The figure also shows the labelling in $T$.
Observe that the finishing labels coincide with the starting ones with
1 and 2 interchanged, so to get the labels for the second twist,
repeat the labels with 1 and 2 interchanged.

The moves from one diagram to the next should all be obvious except
perhaps for the right-most pair in both rows.  Here a ripple spins
round in the surface to eliminate the two oppposite twists (creating
no triple points in the projection). 

Now let $\psi\in\pi_3(BT)$ be the canonical class of this labelled
diagram for $K$.  We shall compute $h(\psi)$ the Hurewicz image in
$H_3(BT)$.  Inspecting the definition of the canonical class given in
section 2 above, we see that $h(\psi)$ is represented by a cycle $C$
given as a sum of 3--cubes of $BT$ one for each triple point of the
diagram.  In figure 3 we have written the 3--cubes (with sign)
determined by the triple points (which appear as R3 moves in the
sequence).  Figure 4 gives more detail of the identification of these
cubes.  Figure 4 shows the cubes $-(012),-(122),(121)$ and $(110)$ the
other two ($-(202)$ and $(102)$) are the same as $-(012)$ and $(121)$
respectively, with label changes.

\figure\epsfxsize4truein
\epsfbox{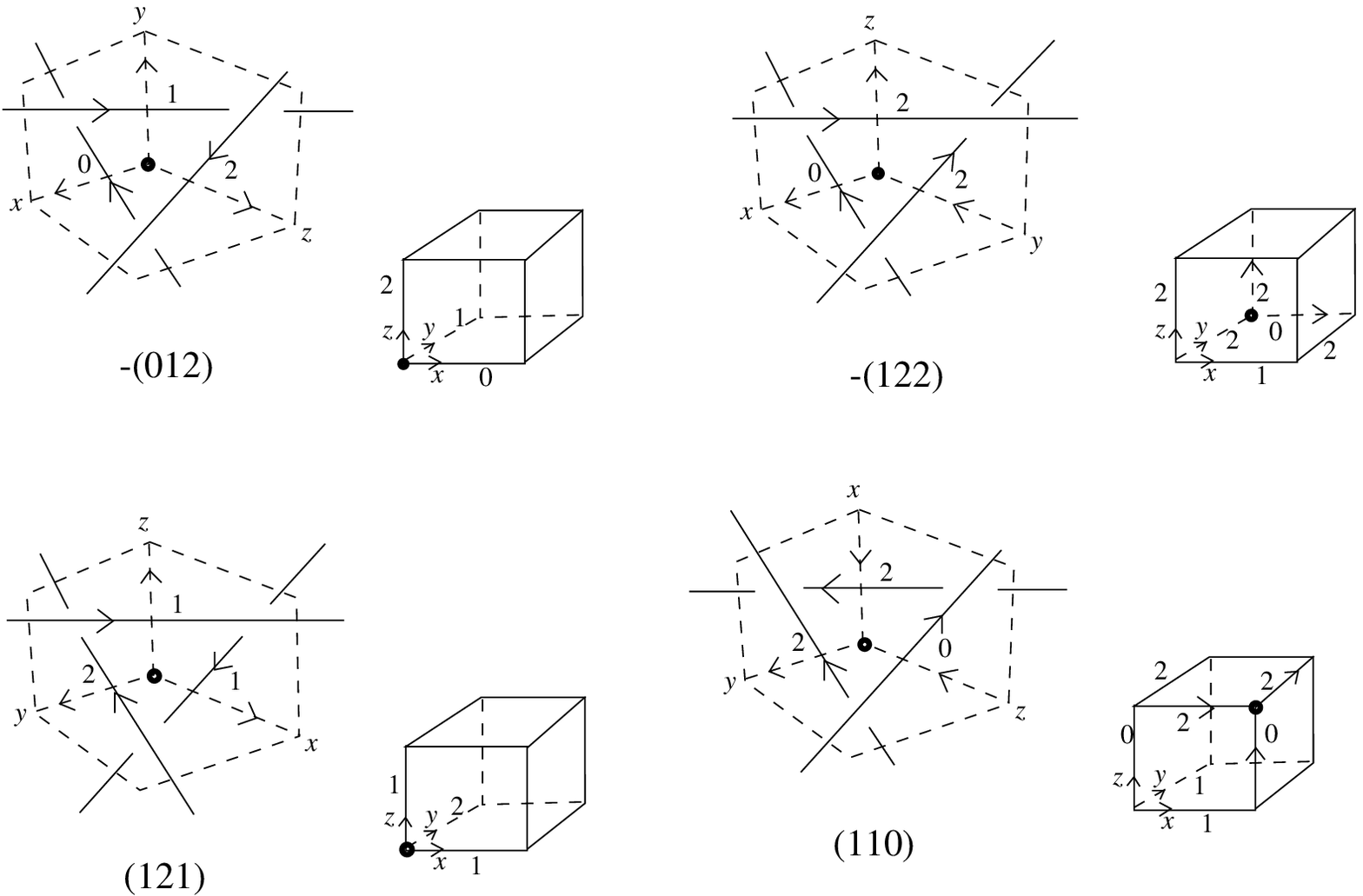}
\endfigure

We can now read off the cycle which represents $h(\psi)$.
$$
\eqalign{C=&-(012)-(202)-(122)+(121)+(102)+(110)\cr
&-(021) -(101) -(211) +(212) +(201) +(220)}
$$
We need two calculations for which details are to
be found in the Maple worksheet [\Maple].

\rk{Calculation 1}$H_3(BT)\cong \Z\times\Z_3$

\rk{Calculation 2}$C$ represents a generator of the $\Z_3$--summand
of $H_3(BT)$.

\proc{Lemma}\key{Invar}The $\Z_3$--summand of $H_3(BT)$ is fixed by
all permutations of $T=\{0,1,2\}$.

\prf Let $Q$ denote the $\Z_3$--summand of $H_3(BT)$ and let
$S_3$ denote the group of permutations of $\{0,1,2\}$.  $C$ represents
a generator of $Q$ and by construction is invariant under the
interchange $(1,2)\in S_3$.  Thus $Q$ is fixed by $(1,2)$.  Now notice
that there are no symmetries of $H_3(BT)$ of order 3. It follows that
the 3--cycle $(0,1,2)\in S_3$ also fixes $Q$.  Since $S_3$ is
generated by $(1,2)$ and $(0,1,2)$, $Q$ is fixed by the whole of
$S_3$.\endprf

\proc{Corollary}\key{Indep}$h(\psi)$ is independent of the choice of 
(non-constant) labelling in $T$. \endproc

\prf  The labelling is clearly determined by the initial labelling of
the first diagram and it follows that two different (non-constant)
labels are related by a permutation of $T$.\endprf

Now let $K'$ denote $K$ with opposite orientation.  Notice that since
we are labelling in an {\sl involutory} rack $T$ (ie one in which
$a^{b^2}=a$ for all $a,b$) the original labels give a labelling after
the orientation change.  Let $\psi'$ be the canonical class of $K'$
using this labelling.  To calculate $h(\psi')$ we observe that
changing orientation reverses all arrows and from figure 5 we see
that this corresponds to replacing each cube $(a,b,c)$ by the cube
$(a^{bc},b^c,c)$ with opposite orientation.  
\figure
\relabelbox\small
\epsfbox{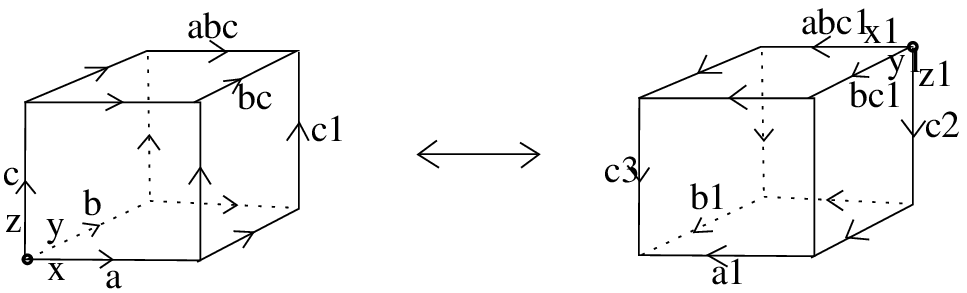}
\relabel {x}{$\ss x$}
\relabel {x1}{$\ss x$}
\relabel {y}{$\ss y$}
\relabel {y1}{$\ss y$}
\relabel {z}{$\ss z$}
\relabel {z1}{$\ss z$}
\relabel {a}{$a$}
\relabel {a1}{$a$}
\relabel {bc}{$b^c$}
\relabel {bc1}{$b^c$}
\relabel {b}{$b$}
\relabel {b1}{$b$}
\relabel {c}{$c$}
\relabel {c1}{$c$}
\relabel {c2}{$c$}
\relabel {c3}{$c$}
\relabel {abc}{$a^{bc}$}
\relabel {abc1}{$a^{bc}$}
\endrelabelbox
\endfigure
Performing this substitution for the cycle $C$ given above we read off
the cycle $C'$ representing $h(\psi')$:
$$
\eqalign{C'=\,\,&(202) +(012)+ (122) -(201) -(212) -(220)\cr
&(101)+ (021)+ (211) -(102) -(121) -(110)}
$$
Serendipitously we notice:

\rk{Observation}$C+C'=0$\endrk

It follows that $h(\psi')=-h(\psi)$.  Now suppose that $K$ is isotopic
to $K'$ preserving orientation.  Then since the fundamental racks of
$K$, $K'$ and the isotopy (thought of an an embedding of $S^2\times I$
in $S^4\times I$) are all isomorphic, the labelling of $K$ in $T$
induces a labelling through the isotopy which is non-trivial on $K'$.
By corollary \ref{Indep} we can use these labels to calculate $h(\psi)$
and $h(\psi')$.  But by lemma \ref{Canon} $\psi'=\psi$ and we have a
contradiction, completing the proof of the main theorem.\endprf

See the remarks in section 5 for explanations of some (but not all)
of the concidences which occured in the above proof.

\section{Trefoils left and right}

We shall need a generalisation of the rack space.  Let $R$ be a rack.
The {\sl extended rack space} $B_RR$ has set of $n$--cubes in
bijection with $R^{n+1}$ (ie the $n+1$--cubes of $BR$) and the same
formulae for face maps as in $BR$ with an index shift (ie,
$\d_i^{\ep}$ in $B_RR$ is $\d_{i+1}^{\ep}$ in $BR$) cf [\Trunks,
1.3.2 and 3.1.2].  It follows that $H_n(B_RR)\cong H_{n+1}(BR)$ (cf
[\James, 5.14 and above]).  Note that $B_RR$ is a covering space of
$BR$.

The interpretation for diagrams is labelling of both arcs and regions.
More precisely, suppose that $D$ is a diagram in $\re^2$ for a knot
$K$ in $\re^3$.  An {\sl extended labelling} of $D$ by a rack $R$ is a
labelling of arcs and regions of $D$ by elements of $R$ with the rules
illustrated in figure 6 for labels of adjacent regions and at
crossings:

\figure
\relabelbox\small
\epsfbox{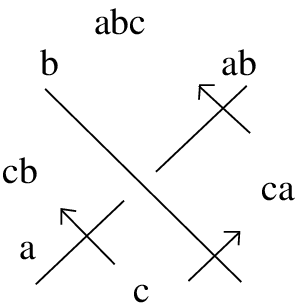}
\relabel {a}{$a$}
\relabel {ab}{$a^b$}
\relabel {b}{$b$}
\relabel {c}{$c$}
\relabel {ca}{$c^b$}
\relabel {cb}{$c^a$}
\relabel {abc}{$c^{ab}=c^{ba^b}$}
\endrelabelbox
\endfigure

Given an extended labelled diagram, there is a map $\re^2\to B_RR$
constructed similarly to section 2: Choose bicollars, map little
squares at crossings to the appropriate 2--cube of $B_RR$ (for example
a square at the crossing in figure 6 would be mapped to the 2--cube
$(c,a,b)$) collar lines across arcs to the appropriate 1--cube (eg, a
collar line across the lower left arc in figure 6 would be mapped to
the 1--cube $(c,a)$) and regions to the 0--cube given by the label.
Thus the labelling determines a canonical class $\phi\in\pi_2(B_RR)$
(based at the vertex corresponding to the label of the infinite
region).  An isotopy of $K$ in $\re^3$ corresponds to a diagram in
$\re^2\times I$ to which the labelling can be canonically extended.
This in turn gives a map $\re^2\times I\to B_RR$, ie a homotopy of the
canonical class.  Thus as in section 2, the canonical class is an
invariant of the isotopy class of the labelled diagram.

We now exhibit extended labellings in the three colour rack $T$ for
diagrams of the right-hand trefoil $K$ and the left-hand trefoil $K'$,
see figure 7.

\figure
\epsfxsize\hsize
\epsfbox{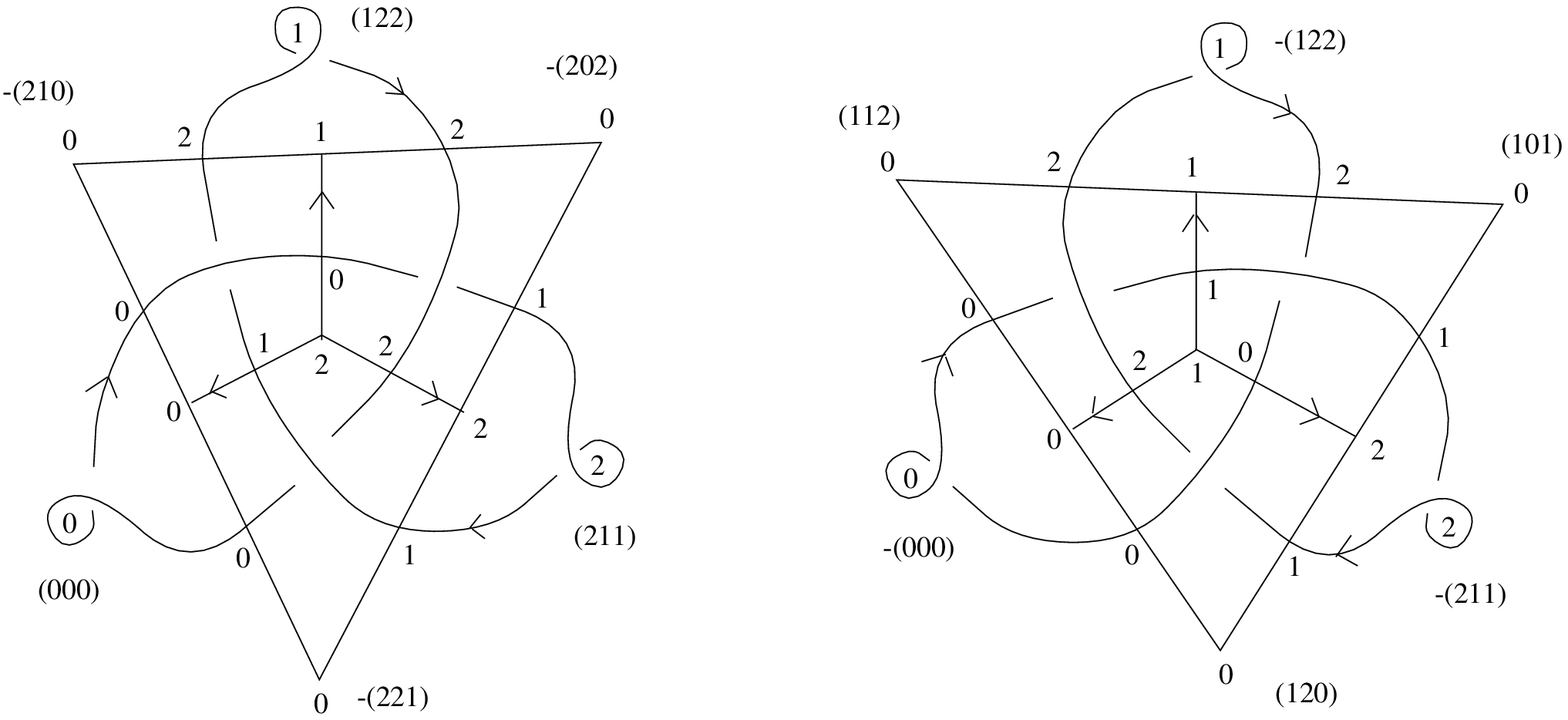}
\endfigure

Let $\phi$ and $\phi'$ be the canonical classes determined by these
diagrams and let $h(\phi)$ and $h(\phi')$ be their Hurewicz images in
$H_2(B_TT)\cong H_3(BT)\cong\Z\times\Z_3$.  We can read off
representating cycles $B$ and $B'$ from the diagrams:
$$
\eqalign{
B=&-(210)-(202)-(221)+(211)+(122)+(000)\cr
B'=&\ (221)+(202)+(210)-(122)-(000)-(211)}
$$
\rk{Observation}$B=-B'$ and hence $h(\phi)=-h(\phi')$.\endrk

We need one final Maple calculation:

\rk{Calculation 3}$B$ represents a generator of the $\Z_3$--summand
of $H_2(B_TT)$.\endrk

\proc{Lemma}\key{TIndep}The classes $h(\phi)$ and $h(\phi')$ are 
independent of the labels in $T$ of arcs and regions in the
diagrams provided that the knots themselves have non-constant
labels.\endproc

\prf We notice that both diagrams have rotational symmetry of order
3.  But any two labellings which are non-constant on the knots are
related by this symmetry followed by a permtuation of $T$.  The result
follows from lemma \ref{Invar}.\endprf

We can now prove that right and left-hand trefoils are different.

\proc{Theorem}There is no isotopy of $\re^3$ which carries the
(unoriented) right-hand trefoil to the left-hand trefoil.\endproc

\prf It is easy to contruct an isotopy carrying an oriented trefoil to 
its orientation reverse.  Thus it suffices to prove that the oriented
left and right trefoils are not isotopic.  The diagrams drawn above
both have zero writhe and if there is an isotopy between the oriented
left and right trefoils then there would be a framed isotopy between
the diagrams.  Now the labelling on the right trefoil determines a
labelling of the isotopy (thought of as a diagram of $S^1\times I$ in
$\re^3\times I$) and hence a labelling of the left trefoil which is
non-constant on the knot.  By lemma \ref{TIndep} we can use these
labels to calculate $h(\phi)$ and $h(\phi')$ and it follows that
$h(\phi)=h(\phi')$ contradicting the observation made above.\endprf

\section{Some remarks}

\sh{Fundamental quandles}

The {\sl involutive} fundamental quandle (ie the fundamental quandle
with the operator relations $a^2\equiv 1$ for each $a$) of the trefoil
is the 3--colour rack/quandle $T$.  This explains why there is a
natural labelling in this rack.  Furthermore using a Van Kampen
argument it can be seen that the fundamental quandle of an
$n$--twist-spun knot $K$ is the fundamental quandle of $K$ with the
operator relations $a^n\equiv 1$ for each $a$.  It follows that the
fundamental quandle of the 2--twist-spun trefoil is $T$.  So the
labelling of the 2--twist-spun trefoil we used is the only non-trivial
labelling (in a quandle).  

\sh{Using the canonical class}

By the Classification Theorem, the canonical class can be used to
distinguish the left and right trefoils which have isomorphic
fundamental racks.  The canonical class lies in $\pi_2(B\Gamma)$ which
is $H_2$ of the universal cover of $B\Gamma$ and therefore determines
classes in every connected cover of $B\Gamma$.  To distinguish the
left and right trefoils we used $B_TT$ which is a 3--fold cover of
$BT$ and the canonical class we used was the image of the class in the
3--fold cover of $B\Gamma$ determined in this way.  This is typical of
the way in which invariants can be constructed to extract information
from the Classification Theorem.  For more examples see [\Mono].

\sh{The $\Z$ factor}

Notice that $H_3(BT)\cong H_2(B_TT)$ has a $\Z$ factor which played no
role in the proofs.  In fact every homology group of the rack space
$BR$ where $R$ has a base element which acts trivially on itself (any
element in a quandle will do) splits a $\Z$ factor which comes from
the inclusion and natural projection $B*\subset BR\to B*$ where $B*$
denotes the rack space of the one element rack.  Now $B*$ is a model
for $\Omega(S^2)$ and has one cell in each dimension with all boundary
maps trivial and hence has homology $\Z$ in each dimension [\James;
3.3].

Now for a 2--knot the canonical class maps to $\pi_3(\Omega
S^2)=\pi_4(S^3)=\Z_2$ which in turn maps to zero in $H_3(\Omega
S^2)=\Z$.  This explains why the $\Z$ factor was immaterial for the
2--twist spun trefoil.  For a 1--knot the image in $H_2(\Omega
S^2)=\Z$ is readily seen to be the writhe of the diagram.  But we
chose our diagrams for the left and right trefoils to have zero
writhe, which explains why the $\Z$ factor was immaterial for this
case as well.

\sh{Invariance under permutations of $T$}

Rack spaces are {\sl simple} in other words $\pi_1$ acts trivially on
$\pi_n$ for each $n$ [\James; proposition 5.2].  The proof is by
diagram manipulation.  Introduce a small sphere labelled by $x\in R$.
Pull the sphere over the diagram (realising the change of labels
corresponding to the action of $x$) and then pull the sphere back
under the diagram and eliminate it.  The same proof shows that any
labelled diagram is bordant to the diagram with labels acted on by
any element of the labelling rack.  Thus the operator group of the
rack acts trivially on the canonical class.  Now the three colour rack
$T$ has operator group equal to its automorphism group, which is the
symmetric group $S_3$.  This explains why the canonical class
determined by a labelling in $T$ is invariant under choice of labels:
any rack with operator group equal to its automorphism group would
have the same property.  Thus the labelling used for the 2--twist-spun
trefoil was immaterial.  For extended labels the same proof shows
invariance under change of labelling with fixed label at infinity.
This together with the obvious symmetry of the choice of label at
infinity in $T$ explains why the labelling was also immaterial for the
left and right trefoils.

\references

\bye